\documentclass[11pt]{article}
\usepackage{amsmath}
\usepackage{amssymb}
\usepackage{amsthm}
\usepackage{amsfonts}
\usepackage{graphicx}
\usepackage{pdfpages}
\usepackage{subcaption}
\usepackage{xcolor,colortbl}
\usepackage{multicol}
\usepackage{url}

\definecolor{Gray}{gray}{0.92}
\definecolor{LightCyan}{rgb}{0.88,1,1}

\newcolumntype{a}{>{\columncolor{Gray}}c}
\newcolumntype{b}{>{\columncolor{white}}c}

\usepackage{bbm,epsfig,graphics,epic,color,rotating,color}
\numberwithin{equation}{section}

\newcommand{\dE}{\mathbb{E}}
\newcommand{\dP}{\mathbb{P}}

\newcommand{\cF}{\mathcal{F}}
\newcommand{\E}{\mathbb{E}}
\newcommand{\V}{\mathbb{V}}

\newtheorem{thm}{Theorem}[section]

\newtheorem{lem}{Lemma}[section]

\newtheorem{cor}[lem]{Corollary}

\usepackage{graphicx}% Include figure files
\usepackage{dcolumn}% Align table columns on decimal point
\usepackage{bm}% bold math
\title{Functional central limit theorem for dependent models with finite memory}
\author{V\'ictor Hugo~V\'azquez-Guevara\thanks{Facultad de Ciencias F\'isico Matem\'aticas, Benem\'erita Universidad Aut\'onoma de Puebla, M\'exico. e-mail: victor.vazquezg@correo.buap.mx} ,
Manuel Gonz\'alez-Navarrete \thanks{Departamento de Matem\'atica y Estad\'istica, Universidad de La Frontera, Chile. e-mail: manuel.gonzaleznavarrete@ufrontera.cl}
\date{}
}

\begin{document}
\maketitle

\begin{abstract}
We provide complementary results for a family of models with dependence on their previous $k$-sum. Using a martingale-based approach, we establish a functional central limit theorem and analyze the limiting behavior of the center of mass. Additionally, we explore the connection between our findings and the study of certain reinforced random walks in the literature.
\end{abstract}


\section{Introduction}

 Before introducing the main model we will deal with, let us consider the generalization of the binomial distribution discussed in \cite{DF}. More precisely, for $n \geq 0$, the random variable $S_n$ counts the number of successes in a correlated Bernoulli sequence, denoted by $\{X_n, n \ge 1\}$, where $X_1$ is distributed as a Bernoulli random variable with parameter $p$ and for $n \geq 1$ correlation is dictated by conditional probabilities defined, for some parameter $0 \leq \theta < 1$, as follows:
\begin{equation}\label{model}
\mathbb{P} ( X_{n+1} = 1 \mid S_n ) = (1-\theta ) p + \theta \frac{S_n}{n}.
\end{equation}
 Of course, we have that $S_n = X_1 + X_2 + \cdots + X_n$ at time $n \geq 1$. The properties and asymptotic analysis of the generalized binomial process expressed in \eqref{model} have been extensively studied in the literature, see \cite{Dre,CBP2024,heyde2004,JJQ,VePu} for details.

We now consider the modified model introduced in \cite{SKV}, where, for a fixed parameter $k \geq 1$, the conditional probability of success in the $i$-th trial depends on the number of successes in all previous trials if $i \leq k + 1$, and on those in the immediate past $k$ trials if $i \geq k + 2$. Specifically, setting $P(X_1 = 1) = p$, we define:
\begin{equation}
\label{k-sum}
\begin{array}{l}
\dP(X_i = 1 \mid \mathcal{F}_{i-1}) = (1 - \theta)p + \frac{\theta}{i - 1} S_{i-1}, \quad \text{for } 2 \leq i \leq k + 1, \\[0.3cm]
\dP(X_i = 1 \mid \mathcal{F}_{i-1}) = (1 - \theta)p + \frac{\theta}{k} S_{i,k}, \quad \text{for } k + 2 \leq i \leq n,
\end{array}
\end{equation}
where $0 \leq \theta < 1$, $0 < p < 1$, $\mathcal{F}_i$ denotes the $\sigma$-field generated by the random variables $X_1, \ldots, X_i$, and $S_{i,k} := S_{i-1} - S_{i-k-1}$. The authors of \cite{SKV} established a central limit theorem, which states that, as $n \to \infty$,
\begin{equation}
\label{CLT}
\frac{S_n - np}{\sqrt{n}} \overset{d}{\longrightarrow} N\left(0, \sigma^2_{p,\theta,k}\right),
\end{equation}
where $\overset{d}{\longrightarrow}$ denotes convergence in distribution and the variance $\sigma^2_{p,\theta,k}$ is given by $\sigma^2_{p,\theta,1} := \frac{p(1 - p)(1 + \theta)}{(1 - \theta)}$ for $k = 1$, and for $k \geq 2$:
\begin{equation}
\label{sigma}
\sigma^2_{p,\theta,k} := \frac{p(1 - p)}{(1 - \theta)^2}
\begin{cases}
\left(1 - \frac{\theta^2}{k^2}
\left\{
 \frac{k - (B(k, 2\theta))^{-1}}{1 - 2\theta}
\right\}\right), & \text{if } \theta \neq \frac{1}{2}, \\[10pt]
\left(1 - \frac{\theta^2}{k} \displaystyle\sum_{j=1}^k \frac{1}{j}\right), & \text{if } \theta = \frac{1}{2},
\end{cases}
\end{equation}
where $B(a,b)$ denotes to the beta function. Additionally, the following law of the iterated logarithm was established:
\[
\limsup_{n \to \infty} \frac{S_n - np}{\sqrt{2n \log \log n}} = \sigma_{p,\theta,k} \quad \text{a.s.},
\]
and
\[
\liminf_{n \to \infty} \frac{S_n - np}{\sqrt{2n \log \log n}} = -\sigma_{p,\theta,k} \quad \text{a.s.}
\]

In this work, we extend the asymptotic analysis for the sequence $(S_n)$ defined by \eqref{k-sum}, demonstrating a proper version of the functional central limit theorem. To achieve this, we employ a novel martingale approach recently developed in \cite{BV2021,GLV2024,CBP2024,gue}. Furthermore, we establish a central limit theorem and the law of large numbers for the center of mass of the number of successes, defined by
\begin{equation}
\label{centerofmass}
C_n = \frac{1}{n} \sum_{i=1}^n S_i,
\end{equation}
and introduced in \cite{Grill,LoWa}. In addition, we draw connections between the dependency structure expressed in \eqref{k-sum} and other stochastic processes such as the elephant random walk \cite{Schutz} and the minimal random walk \cite{kumar}.

The rest of the paper is structured as follows. The next section establishes a connection between the considered dependent model and reinforced random walks studied in recent literature. The main results are presented in Section \ref{sec:main}. Finally, Section \ref{proofs} is dedicated to the proofs of our findings.

\section{Dependent models and reinforced random walks}
\label{sec:RWs}

In this section we discuss the relation of proposed model with some reinforced random walks. It is worth mentioning that there exists an emergent literature on such reinforced random walks, including  problems in different contexts, see for instance \cite{BercuHyp,BL,Bertoin3,GLR}.

\subsection{The minimal random walk}
\label{sec:MRW}
First, denote by \( (X^M_n) \) a random walk which at $i$-th step chooses a random time $U(i)$ from its past and perform the $(i+1)$-th step by following the distributional rule
\begin{equation}
\label{MRW}
X^M_{i+1} \sim 
\begin{cases}
Ber(r), & \text{if } X_{U(i)} = 1, \\[10pt]
Ber(q), & \text{if } X_{U(i)} = 0,
\end{cases}
\end{equation}
that is, a Bernoulli random variable with parameters \( r \) and \( q \) in each case. The case in which \( U(i) \) is uniformly distributed from the whole past is well known, see for instance \cite{BV2021,VCG}. We define an evolution inspired by the conditional probabilities in \eqref{k-sum}. That is, in the case \( i \leq k +1 \), \( U(i) \sim Unif(\{1, 2, \ldots, i\}). \) Otherwise, if \( i > k +1 \), the walker considers the previous \( k \) steps, that is \( U(i) \sim Unif(\{i-k+1, i-k+2, \ldots, i\}). \)

This formulation is inspired by the so-called minimal random walk \cite{kumar}. In this sense, its conditional probabilities are given by
\begin{equation}
\label{dPMRW}
\begin{array}{l}
\dP(X_i^M = 1 \mid \mathcal{F}_{i-1}) = r\frac{ S_{i-1}^M}{i-1} + q\left(\frac{ i-1-S_{i-1}^M}{i-1}\right) = q + (r-q)\frac{S_{i-1}^M}{i-1}, \quad \text{for } 2 \leq i \leq k + 1, \\[0.3cm]
\dP(X_{i}^M =1 \mid \mathcal{F}_{i-1}) = r\frac{ S_{i,k}^M}{k} + q\left(\frac{ k-S_{i,k}^M}{k}\right) = q + (r-q)\frac{S_{i,k}^M}{k}, \quad \text{for } k + 2 \leq i \leq n.
\end{array}
\end{equation}
Then, the variance of the process is given by \( \sigma^2_{r,q,1} = \frac{q(1-r)(1+r-q)}{(1-(r-q))^3}, \) and for \( k \geq 2 \),
\[
\sigma_{r,q,k}^{2} = \frac{q(1-r)}{(1-(r-q))^4}
\left\{\begin{array}{ccc}\left(1-\frac{\left(r-q\right)^2}{k^2}\left(\frac{k-[B(k,2(r-q))]^{-1}}{1-2(r-q)}\right)\right)& , & \text{ if }  r-q \neq \frac{1}{2},  \\ \left(1-\frac{\left(r-q\right)^2}{k}\displaystyle\sum_{j=1}^k \frac{1}{j}\right)&,& \text{ if } r-q = \frac{1}{2}.\end{array}
\right.
\]

\subsection{The elephant random walk}
\label{sec:ERW}

On the same fashion, we define a variation of the elephant random walk, which at $i$-th step draws a random variable $X_i^{E} \in \{-1,+1\}$ by choosing a step from its past in the same form than the finite memory minimal random walk, however its current step follows the rule:
\begin{equation}
\label{ERW}
X^E_{i+1} = 
\begin{cases}
X^E_{U(i)}, & \text{with probability} \ \alpha, \\[10pt]
-X^E_{U(i)}, & \text{with probability} \ 1-\alpha.
\end{cases}
\end{equation}
Let denote $S_n^{E}$ the position of the elephant, it is possible to observe that its conditional probabilities are 
\begin{equation}
\label{dPERW}
\begin{array}{l}
\dP(X_i^E = 1 \mid \mathcal{F}_{i-1}) =  \frac{1}{2}\left( 1+  (2\alpha-1)\frac{S_{i-1}^E}{i-1} \right), \quad \text{for } 2 \leq i \leq k + 1, \\[0.3cm]
\dP(X_{i}^E =1 \mid \mathcal{F}_{i-1}) =  \frac{1}{2}\left( 1+  (2\alpha-1)\frac{S_{i,k}^E}{k} \right), \quad \text{for } k + 2 \leq i \leq n.
\end{array}
\end{equation}
In this case, one may observe that the variance of the process is given by \( \sigma^2_{\alpha,1} = \frac{\alpha}{(1-\alpha)} \), and for \( k \geq 2 \),
\[
\sigma_{\alpha,k}^{2} = \frac{1}{4(1-\alpha)}
\left\{\begin{array}{ccc}\left(1-\frac{\left(2\alpha-1\right)^2}{k^2}\left(\frac{k-[B(k,2(2\alpha-1))]^{-1}}{3-4\alpha}\right)\right)& , & \text{ if }  \alpha \neq \frac{3}{4},  \\ \left(1-\frac{\left(2\alpha-1\right)^2}{k}\displaystyle\sum_{j=1}^k \frac{1}{j}\right)&,& \text{ if } \alpha= \frac{3}{4}.\end{array}
\right.
\]

This model is related to a variation of the elephant random walk discussed in \cite{Gut}. In such paper, the authors showed asymptotic normality as in \eqref{CLT} for the cases $k=1$ (Theorem 7.1) and $k=2$ (Theorem 8.1), by using the Markov chain approach in \cite{Jo}. In next section we complement their results.

\section{Main Results}
\label{sec:main}

In this section we present the main results of this work for the dependent model with finite memory. The following functional central limit theorem complements Theorem 3.5 of \cite{SKV}.

\begin{thm}
\label{T_Funct}
We have the distributional convergence in $D([0,\infty[)$ the Skorokhod space of right-continuous functions with left-hand limits, as $n \to \infty$,
\begin{equation}
\label{FCLT-DR}
\left( \sqrt{n}\Big(\frac{S_{\lfloor nt \rfloor}}{\lfloor nt \rfloor}-p\Big), t \geq 0\right) \overset{d}{\longrightarrow} \big( W_t, t \geq 0 \big),
\end{equation}
where $\big( W_t, t \geq 0 \big)$ is a real-valued centered Gaussian process starting at the origin with covariance given by $\dE[W_s W_t]= \frac{1}{t}\sigma^2_{p,\theta,k}$, for all $0<s \leq t$.
In particular, as $n \to \infty$,
\[\sqrt{n}\Big(\frac{S_n}{n}-p\Big) \overset{d}{\longrightarrow} N\left(0,\sigma^2_{p,\theta,k}\right).\]
\end{thm}

Furthermore, we also provide the law of large numbers and the central limit theorem for the center of mass of the number of successes, as defined in \eqref{centerofmass}.

\begin{thm}
\label{CLTCM}
We have the almost sure convergence,
\begin{equation}
\label{SLLN-CM-DR}
\lim_{n \rightarrow \infty} \frac{C_n}{n}=\frac{p}{2},
\end{equation}
and the asymptotic normality, as $n \to \infty$,
\begin{equation}
\label{CLT-CM}
\sqrt{n}\left(\frac{C_n}{n}-\frac{p}{2}\right) \overset{d}{\longrightarrow} N\left(0,\frac{\sigma^2_{p,\theta,k}}{3}\right).
\end{equation}
\end{thm}

 In what follows we enunciate the relation of this model with the random walks in Section \ref{sec:RWs}.
 
 In the case of the finite memory minimal random walk, we get the following results.

 \begin{cor}
\label{cor:MRW}
We have, as $n \to \infty$,
\[\left\{\sqrt{n}\left(\frac{S_{\lfloor nt\rfloor}^M}{\lfloor nt\rfloor} - \frac{q}{1-r+q}\right) \ , \ t\geq 0 \right\} \overset{d}{\longrightarrow} \left\{W_t \ , \ t\geq 0\right\},\]
where $\E[W_sW_t] = \frac{1}{t}\sigma^2_{r,q,k}$ for all $0<s\leq t$.

Moreover, as $n \to \infty$, we have that its center of mass $C^M_n \to \frac{q}{2(1-r+q)}$ almost surely and
\begin{equation}
\sqrt{n}\left(\frac{C^M_n}{n}-\frac{q}{2(1-r+q)}\right) \overset{d}{\longrightarrow} N\left(0,\frac{\sigma^2_{r,q,k}}{3}\right).
\end{equation}
The details of $\sigma_{r,q,k}^2$ are given in Section \ref{sec:MRW}.
 \end{cor}
 
 In the case of the finite memory elephant random walk, we get the following results.
 \begin{cor}
\label{cor:ERW}
We have, as $n \to \infty$,
\[\left\{\sqrt{n}\left(\frac{S_{\lfloor nt\rfloor}^E}{\lfloor nt\rfloor} \right) \ , \ t\geq 0 \right\} \overset{d}{\longrightarrow} \left\{W_t \ , \ t\geq 0\right\},\]
where $\E[W_sW_t] = \frac{1}{t}\sigma^2_{\alpha,k}$ for all $0<s\leq t$.

Moreover, as $n \to \infty$, we have that almost surely: $C^E_n \to 0$  and
\begin{equation}
\frac{C^E_n}{\sqrt{n}}\overset{d}{\longrightarrow} N\Big(0,\frac{\sigma^2_{\alpha,k}}{3}\Big).
\end{equation}
The details of $\sigma_{\alpha,k}^2$ are given in Section \ref{sec:ERW}.
 \end{cor}

\section{Proofs}\label{proofs}

The proofs are based on discrete time martingale theory, hence we define a martingale $(M_n)$ inspired by the one considered in \cite{SKV} and we will characterize its behaviour by the martingale differences $(L_n)$, via the relation
\begin{equation} \label{martingale2}
M_n=\sum_{i=1}^n L_i.
\end{equation}
In the original paper \cite{SKV} it was taken into account the martingale difference sequence $\{L_n, n \geq 1\}$, given by $L_1 = X_1 - p$ and 
\[
L_i =
\begin{cases} 
X_i - (1 - \theta)p - \frac{\theta}{i - 1} S_{i-1}, & \text{for } 2 \leq i \leq k + 1, \\[10pt]
X_i - (1 - \theta)p - \frac{\theta}{k} S_{i,k}, & \text{for } k + 2 \leq i \leq n.
\end{cases}
\]
Furthermore, note that almost surely 
\[L_{n+1}=X_{n+1}-\E[X_{n+1}|\mathcal{F}_n],\]
and, since $X_{n+1} \in \{0,1\}$, then $\mathbb{E}[X_{n+1}^l|\mathcal{F}_{n}] = \mathbb{E}[X_{n+1}|\mathcal{F}_{n}]$, for all $l \ge 1$. Additionally, for $l \geq 2$, we find that
\begin{equation}
\mathbb{E}[L_{n+1}^l \mid \mathcal{F}_n] = \sum_{j=0}^{l-2} {l \choose j} (-1)^j \left(\mathbb{E}[X_{n+1}|\mathcal{F}_{n}]\right)^{j+1} 
+ (-1)^{l-1}(l-1)\left(\mathbb{E}[X_{n+1}|\mathcal{F}_{n}]\right)^l 
\hspace{1cm} \text{a.s.,}
\notag
\end{equation}
which implies that
\begin{equation} \label{momentosl}
\sup_{n \geq 0} \mathbb{E}\left[L^l_{n+1} \mid \mathcal{F}_n\right] < \infty, \ \text{a.s.,}
\end{equation}
because $\E[X_{n+1} | \mathcal{F}_n] \leq 1$, almost surely. In particular, 
\begin{equation} \label{segunda}
\mathbb{E}\left[L^2_{n+1}|\mathcal{F}_n\right]=\left(\mathbb{E}[X_{n+1}|\mathcal{F}_{n}]\right)\left(1-\left(\mathbb{E}[X_{n+1}|\mathcal{F}_{n}]\right)\right),
\ \text{a.s.}
\end{equation}
Hence, let denote $e_n := \E[X_{n+1}|\mathcal{F}_{n}]$, then by analysing the polynomial $p(e_n)= e_n-e_n^2$ we have that
\begin{equation} \label{cotasegundomom}
\sup_{n\geq 0} \mathbb{E}\left[L^2_{n+1}|\mathcal{F}_n\right] \leq \frac14,  \ \text{a.s.}
\end{equation}
On the same line, it may be seen that almost surely
\begin{equation}
\label{M4EPS}
\dE[L_{n+1}^4|\cF_n]= e_n-4e_n^2+6e_n^3-3e_n^4.
\notag
\end{equation}
By maximizing this polynomial, we conclude that
\begin{equation}\label{cuarta}
\sup_{n\geq 0} \mathbb{E}\left[L^4_{n+1}|\mathcal{F}_n\right] \leq \frac1{12}, \ \text{a.s.}
\end{equation}

\subsection{Proof of Theorem \ref{T_Funct}}

In order to demonstrate the functional central limit theorem by applying Theorem 1 of \cite{Touati} or Theorem A.9 of \cite{ThesisLucile}, we must analyze the predictable quadratic variation of {$(M_n)$} defined by $\langle M \rangle_0=0$ and for all $n \geq 1$, by
\begin{equation} \label{procrec}
\langle M \rangle_n = \sum_{i=1}^n \mathbb{E}[L_i^2| \mathcal{F}_{i-1}].
 \end{equation}
In addition, during the proof of Theorem 3.5 of \cite{SKV} it was achieved that
\[
\frac{\langle M\rangle_n}{n} \to (1-\theta)^2 \sigma^2_{p,\theta,k} \ \ \mbox{in probability,}\]
which leads to 
\[
\frac{1}{n}\langle M\rangle_{\lfloor nt\rfloor} \to t (1-\theta)^2 \sigma^2_{p,\theta,k}  \ \ \mbox{in probability.}
\]
In order to prove the corresponding Lindeberg's condition, we have for any $\varepsilon > 0$ that 
\[
\displaystyle\frac{1}{n}\sum_{i=1}^n\dE [L_i ^2 \mathbb{I}_{\{|L_i| > \varepsilon \sqrt{n}\}} \vert \mathcal{F}_{i-1}] \leq  \displaystyle\frac{1}{n^2\varepsilon^2}\sum_{i=1}^n\dE [L_i ^4 \vert \mathcal{F}_{i-1}] \leq  \displaystyle\frac{1}{12 n^2\varepsilon^2}n,
\]
where first inequality holds since we are in the event $\frac{|L_i|^2}{(\varepsilon \sqrt{n})^2} \geq 1$, while second inequality follows from \eqref{cuarta}. Therefore, we obtain that
\[
\displaystyle\frac{1}{n}\sum_{i=1}^n\dE [L_i ^2 \mathbb{I}_{\{|L_i| > \varepsilon \sqrt{n}\}} \vert \mathcal{F}_{i-1}] \rightarrow 0 \text{ as } n \to \infty \text{ in probability}.
\]
Then, we conclude that for all $t\ge 0$ and for any $\varepsilon > 0$, 
\begin{equation}
\label{Lindfunctional}
\displaystyle\frac{1}{n}\sum_{i=1}^{\lfloor nt \rfloor}\dE [L_i ^2 \mathbb{I}_{\{|L_i| > \varepsilon \sqrt{n}\}} \vert \mathcal{F}_{i-1}] \rightarrow 0,
\end{equation}
as $n \to \infty$ in probability. Hence, from Theorem 1 of \cite{Touati} or Theorem A.9 of \cite{ThesisLucile} we have that
\[\left( \frac{M_{\lfloor nt \rfloor}}{\sqrt{n}}, t\geq0 \right) \overset{d}{\longrightarrow} \left( B_t,  t\geq0 \right),\]
where $\left( B_t,  t\geq0 \right)$ is a real-valued centered Gaussian process starting at the origin with covariance given, for all $0<s\leq t$ by
\[\mathbb{E}\left[B_s B_t\right]=s(1-\theta)^2 \sigma^2_{p,\theta,k}.\]
In addition, equation (20) in \cite{SKV}, implies that
\begin{equation}
\label{Mntfunctional}
 \frac{M_{\lfloor nt \rfloor}}{\sqrt{n}}  \sim \frac{(1-\theta){\lfloor nt \rfloor}}{\sqrt{n}} \left( \frac{S_{\lfloor nt \rfloor}}{{\lfloor nt \rfloor}} -p \right) \hspace{1cm}\text{a.s.}
\end{equation}
This leads us to the weak convergence
\[\left( \sqrt{n} \left( \frac{S_{\lfloor nt \rfloor}}{{\lfloor nt \rfloor}} -p\right), t \geq 0\right) \overset{d}{\longrightarrow} \big( W_t, t \geq 0 \big),\]
where $W_t = B_t /(t(1-\theta))$, is such that for $0<s\leq t$: $\mathbb{E}\left[W_s W_t \right]=\frac{\sigma^2_{p,\theta,k}}{t}$, which completes the proof.

\subsection{Proof of Theorem \ref{CLTCM}}

In order to prove convergence \eqref{SLLN-CM-DR} we use the strategy proposed in \cite{BV2021} by employing the law of large numbers established in Theorem 3.2 of \cite{SKV} that states that, almost surely
\[\frac{S_n}{n}\rightarrow p.\]
Hence, the Toeplitz lemma lets us conclude that
\[\frac{\displaystyle\sum_{i=1}^n i \frac{S_i}{i}}{\displaystyle\sum_{i=1}^n i}\rightarrow p \hspace{1cm} \text{a.s.}\]
Therefore
\[\frac{C_n}{n}=\frac{1}{n^2}\sum_{i=1}^n S_i\rightarrow \frac{p}{2} \hspace{1cm} \text{a.s.}\]
For proving \eqref{CLT-CM}, first of all, note that 
\[
C_n=\int_0^1 S_{{\lfloor nt \rfloor}}dt,
\]
and that $\frac{t}{\lfloor nt \rfloor}\sim \frac{1}{n},$ hence
\[\frac{C_n}{n}=\int_0^1 \frac{S_{{\lfloor nt \rfloor}}}{n}dt \sim \int_0^1 t \frac{S_{{\lfloor nt \rfloor}}}{\lfloor nt \rfloor}dt,\]
which implies that
\[\sqrt{n}\left(\frac{C_n}{n}-\frac{p}{2}\right) \sim \int_0^1 \sqrt{n} \left(\frac{S_{{\lfloor nt \rfloor}}}{\lfloor nt \rfloor}-p\right)t dt, \]
that, thanks to Theorem \ref{T_Funct} converges weakly to \(\displaystyle\int_0^1 t W_t dt,\) whose distribution is $N(0,\nabla^2)$, where
\[\nabla^2 = \dE\Big[\Big( \int_0^1 t W_t dt\Big)^2\Big]=2 \int_0^1 \int_0^t st \dE[W_s W_t] ds dt = 2\sigma^2_{p,\theta,k} \int_0^1 \int_0^t s ds dt=\frac{\sigma^2_{p,\theta,k} }{3}.\]

\subsection{Proof of Corollary \ref{cor:MRW}}
It is enough to note that its conditional probabilities are analogous to \eqref{k-sum}, with the following parameter relation, \( \theta=r-q \) and \( p=\frac{q}{1-(r-q)}\).

\subsection{Proof of Corollary \ref{cor:ERW}}
The proof is straightforward by noticing that the conditional probabilities \eqref{dPERW} are similar to \eqref{k-sum} with $\theta = 2\alpha-1$ and $p=1/2$. Moreover $\E(S_n^{E}) = 2 \E(S_n) -n$ and $\V(S_n^{E}) = 4 \V(S_n)$.


\medskip

\begin{thebibliography}{99}
\bibitem{DF}
Z. Drezner and N. Farnum, “A generalized binomial distribution,” Comm. Statist. Theory Methods \textbf{22}, 3051–3063 (1993).

\bibitem{Dre}
Z. Drezner, “On the Limit of the Generalized Binomial Distribution,” Comm. Statist. Theory Methods \textbf{35}(2), 209–221 (2006).

\bibitem{CBP2024}
M. González-Navarrete, R. Lambert, and V. H. Vázquez-Guevara, “A complete characterization of a correlated Bernoulli process,” Electron. Commun. Probab. \textbf{29}(69), 1–12 (2024).

\bibitem{heyde2004}
C. C. Heyde, “Asymptotics and criticality for a correlated Bernoulli process,” Aust. N. Z. J. Statist. \textbf{46}, 53–57 (2004).

\bibitem{JJQ}
B. James, K. James, and Y. Qi, “Limit theorems for correlated Bernoulli random variables,” Stat. Probab. Lett. \textbf{78}, 2339–2345 (2008).

\bibitem{VePu}
P. Vellaisamy and A. P. Punnen, “On the nature of the binomial distribution,” J. Appl. Probab. \textbf{38}(1), 36–44 (2001).

\bibitem{SKV}
D. Singh, S. Kumar, and P. Vellaisamy, “The limit theorems for a previous k-sum dependent model,” J. Math. Anal. Appl. \textbf{487}(2), 124004 (2020).

\bibitem{BV2021}
B. Bercu and V. H. Vázquez Guevara, “Further results on the minimal random walk,” J. Phys. A: Math. Theor. \textbf{55}(41), 415001 (2022).

\bibitem{GLV2024}
M. González-Navarrete, R. Lambert, and V. H. Vázquez-Guevara, “On the asymptotic analysis of lazy reinforced random walks: a martingale approach,” J. Math. Anal. Appl. \textbf{549}(2), 129520 (2025).

\bibitem{gue}
V. H. Vázquez Guevara, “On the almost sure central limit theorem for the elephant random walk,” J. Phys. A: Math. Theor. \textbf{52}(1), 475201 (2019).

\bibitem{Grill}
K. Grill, “On the average of a random walk,” Statist. Probab. Lett. \textbf{6}(5), 357–361 (1988).

\bibitem{LoWa}
C. H. Lo and A. R. Wade, “On the centre of mass of a random walk,” Stochastic Process. Appl. \textbf{129}(11), 4663–4686 (2019).

\bibitem{Schutz}
G. M. Schütz and S. Trimper, “Elephants can always remember: Exact long-range memory effects in a non-Markovian random walk,” Phys. Rev. E \textbf{70}, 045101 (2004).

\bibitem{kumar}
U. Harbola, N. Kumar, and K. Lindenberg, “Memory-induced anomalous dynamics in a minimal random walk model,” Phys. Rev. E \textbf{90}, 022136 (2014).

\bibitem{BercuHyp}
B. Bercu, M. Chabanol, and J.-J. Ruch, “Hypergeometric identities arising from the elephant random walk,” J. Math. Anal. Appl. \textbf{480}(1), 123360 (2019).

\bibitem{BL}
B. Bercu and L. Laulin, “On the center of mass of the elephant random walk,” Stochastic Process. Appl. \textbf{133}, 111–128 (2021).

\bibitem{Bertoin3}
J. Bertoin, “Counting the zeros of an elephant random walk,” Trans. Amer. Math. Soc. \textbf{375}, 5539–5560 (2022).

\bibitem{GLR}
H. Guérin, L. Laulin, and K. Raschel, “A fixed-point equation approach for the superdiffusive elephant random walk,” Ann. Inst. Henri Poincaré (B) Probab. Stat., to appear.

\bibitem{VCG}
V. H. Vázquez Guevara, C. F. Coletti, and M. González-Navarrete, “Complementary asymptotic analysis for a minimal random walk,” J. Math. Phys. \textbf{66}(5), 053301 (2025).

\bibitem{Gut}
A. Gut and U. Stadtmüller, “Variations of the Elephant Random Walk,” J. Appl. Probab. \textbf{58}(3), 805–829 (2021).

\bibitem{Jo}
G. L. Jones, “On the Markov chain central limit theorem,” Prob. Surveys \textbf{1}, 299–320 (2004).

\bibitem{Touati}
A. Touati, “Sur la convergence en loi fonctionnelle de suites de semimartingales vers un mélange de mouvements browniens,” Teor. Veroyatnost. i Primenen. \textbf{36}(4), 744–763 (1991).

\bibitem{ThesisLucile}
L. Laulin, About the elephant random walk, Ph.D. thesis, Université de Bordeaux (2022), Numerical Analysis [cs.NA].



\end{thebibliography}
\end{document}